\documentclass[12pt,twoside]{report}
\usepackage{amsmath}
\usepackage{amssymb}
\usepackage{amscd}
\usepackage{cases}
\usepackage{enumerate}

\title{Finite-time stability for differential inclusions with applications to neural networks}
\author{Rados\l{}aw Matusik, Andrzej Nowakowski, S\l{}awomir Plaskacz, Andrzej Rogowski}
\newtheorem{thm}{Theorem}[chapter]
\newtheorem{prop}[thm]{Proposition}
\newtheorem{rem}[thm]{Remark}
\newtheorem{cor}[thm]{Corollary}

\newenvironment{prf}{\textbf{Proof.}}{\hfill$\Box$\\}
\usepackage{geometry}
\date{}
\begin{document}
\maketitle
\begin{abstract}
The paper investigates sufficient conditions on a differential inclusion which guarantee that the origin is a finite time stable equilibrium, namely a weak local one, a weak global one or a strong local one. The analysis relies on the existence of a Lyapunov function. A new Gronwall type results are used to estimate the settling time. An example of a neural network which is finite-time stable is given.\\
\textbf{Keywors:} finite-time stability, differential inclusions, neural networks, viability, Gronwall's inequality, contingent derivative\\
\textbf{e-mail:} radoslaw.matusik@wmii.uni.lodz.pl, andrzej.nowakowski@wmii.uni.lodz.pl, plaskacz@mat.umk.pl, andrzej.rogowski@wmii.uni.lodz.pl
\end{abstract}
\chapter{Introduction}
In the paper we consider the problem of finite-time stability and stabilization for nonautonomous differential inclusions $x'\in F(t,x)$. We assume that the set-valued map $F$ is merely $\mathcal{L}\times\mathcal{B}$ - measurable and upper semicontinuous in the state variable $x$. The finite-time stability problem is considered in the weak and in the strong framework. The equilibrium point is strong finite-time stable if every solution reaches the equilibrium in finite time. The equilibrium is weak finite-time stable if for every initial condition there exists a solution reaching the equilibrium in finite time. The weak finite-time stability is related to the problem of stabilization of a control system in finite time. In the paper we present sufficient conditions to weak and strong finite-time stability using Lyapunov second method. The sufficient condition to weak finite-time stability has the form
\begin{equation*}
\inf\limits_{f \in F(t,x)} D_{\uparrow} V(t,x)(1,f)\leq -c(t)g(V(t,x)),
\end{equation*}
where $V(t,x)$ is a Lyapunov function and $D_{\uparrow} V$ denotes its contingent epiderivative. Strong finite-time stability is obtained under the assumption
\begin{equation*}
\sup\limits_{f \in F(t,x)} D_{\downarrow} V(t,x)(1,f)\leq -c(t)g(V(t,x)),
\end{equation*}
where $D_{\downarrow} V$ denotes the contingent hypoderivative (the definition and properties of nonsmooth analysis tools used in the paper can be find in \cite{AF} and \cite{FPR}). We assume that the function $c(\cdot)$ is locally integrable and $\int_{t_0}^\infty c(t)dt>0$ for every $t_0 \geq 0$, the function $g(\cdot)$ is increasing, $g(0)=0$ and $\int_0^b \frac{1}{g(v)}dv<\infty$ for $b>0$. To obtain a global stability result we assume that the function $V$ is radially unbounded. To obtain local stability we merely need that $V$ is positively definite. When we consider weak stability, we assume that the function $V(t,x)$ is epigraphically absolutely continuous in time and continuous in $(t,x)$. This assumption is prompted by our method of proof that bases on Theorem 4.4 in \cite{FP} that has been obtained as a corollary of a viability theorem. This regularity assumption on $V$ is weaker then the assumption on $V$ that we need to obtain strong stability. To receive strong stability we assume that $V$ is locally Lipschitz continuous.

We apply the obtained finite-time stability criteria to show finite-time stability for a class of Hopfield neural networks with discontinuous activation functions.

The problem of finite-time stability of differential inclusions has been considered in the strong framework for autonomous differential inclusions in \cite{MPa}. In \cite{BB} the problem of finite-time stability has been considered for autonomous and continuous differential equation. In \cite{MP} sufficient conditions of finite-time stability has been obtained for nonautonomous differential equation.
\chapter{Finite-time stability of a measurable upper semicontinuous differential inclusion}
\setcounter{equation}{0}
Our aim in the section is to obtain a sufficient condition for weak and strong finite-time stability for a differential inclusion in a local and global framework. We recall Theorem 4.4 from \cite{FP} that concerns the monotone behavior of a solution to a differential inclusion with respect to a continuous Lyapunov function. Theorem 4.4 is a consequence of viability results presented in \cite{FPR} and \cite{FP}. By using a modified version of the Gronwall inequality we obtain the weak finite-time stability for a measurable upper semicontinuous differential inclusion.

We consider a differential inclusion
\begin{equation}\label{2.01}
x'(t)\in F(t,x(t)).
\end{equation}

We shall assume that the origin is an equilibrium point, i.e. $0\in F(t,0)$ for a.a. $t\in[0,\infty)$.

If $x:[t_0,T)\to \mathbb{R}^{n}\setminus\{0\}$ is a solution to (\ref{2.01}) and $\lim_{t\to T^-} x(t)=0$, then $T$ is the settling time of the solution.

If the right-hand side $F$ is single valued and the differential equation has a unique solution in forward time starting at $t_0$ from initial point $x_0$, different to the origin then the settling time is a function of $(t_0,x_0)$. In the autonomous case we set by default $t_0=0$ and the settling time is a function of the initial state $x_0$. The settling-time function has surprising properties even for autonomous differential equation. In \cite[p. 756]{BB} it is provided an example of an autonomous vector field on the plane such that there exists a sequence of initial points converging to the origin for which the settling times tend to infinity.

We say that the origin is a weak local finite-time stable equilibrium point of the differential inclusion (\ref{2.01}) if for every $\varepsilon>0$ and every $t_0\geq 0$ there exists $\delta>0$ such that for every $x_0$ in $\delta B$ there exists an absolutely continuous solution $x:[t_0,\infty)\to\varepsilon B$ of (\ref{2.01}) satisfying $x(t_0)=x_0$ such that $x(t)=0$ for $t\geq T(t_0,x_0)$, where $B$ denotes the unit ball in $\mathbb{R}^n$. It is obvious that $T(t_0,x_0)$ is an upper estimation of the settling time of the solution $x(\cdot)$.

We say that the origin is a weak global finite-time stable equilibrium point if for every initial condition ($t_0,x_0)$ there exists a solution $x(\cdot)$ reaching the origin in finite time.

We say that the origin is a strong local finite-time stable equilibrium point of the differential inclusion (\ref{2.01}) if for every $\varepsilon>0$ and every $t_0\geq 0$ there exists $\delta>0$ such that for every $x_0$ in $\delta B$, every absolutely continuous solution $x:[t_0,\infty) \to \mathbb{R}^n$ of (\ref{2.01}) satisfying $x(t_0)=x_0$ is bounded by $\varepsilon$ and $x(t)=0$ for $t\geq T(t_0,x_0)$. The estimation $T(t_0,x_0)$ of the settling time is common for all solutions satisfying the initial condition $x(t_0)=x_0$.

We assume that the right-hand side $F:[0,\infty)\times \mathbb{R}^{n}\multimap \mathbb{R}^{n}$ satisfies the following conditions:
\begin{equation}\label{2.02}
\left\{
\begin{array}{l}
F(t,x) \textup{ is a nonempty compact convex set for } (t,x) \in [0,\infty)\times \mathbb{R}^{n}, \\
F(t,\cdot) \textup{ is upper semicontinuous for almost all t}, \\
F(\cdot,\cdot) \textup{ is } {\mathcal L}(\mathbb{R}) \times {\mathcal B}(\mathbb{R}^{n}) \textup{ (Lebesgue-Borel) - measurable}, \\
|F(t,x)| \leq \mu(t)(1+|x|),
\end{array}
\right.
\end{equation}
where $\mu(\cdot)$ is a locally integrable function and $|F(t,x)|=\sup\{|f|:f\in F(t,x)\}$. Recall that if $F$ is a Carath\'{e}odory map in the sense that for every $x\in \mathbb{R}^{n}$, $F(\cdot,x)$ is measurable and for almost every $t\in[0,\,\infty)$, $F(t,\cdot)$ is continuous, then it is Lebesgue-Borel measurable. We say that a tube (a set-valued map) $P:[0,\infty]\multimap \mathbb{R}^{n}$ is absolutely continuous if for every compact set $K\subset \mathbb{R}^{n}$ there exists a locally integrable function $\mu_K:[0,\infty)\to[0,\infty)$ such that for every $s \in [0,t)$ we have
\begin{gather}\label{2.03}
\max\{e(P(t)\cap K,P(s)),e(P(s)\cap K,P(t))\} \leq\int\limits_s^t \mu_K(\tau)d\tau,
\end{gather}
where the Hausdorff semi-distance of sets $A$, $C$ is given by $e(A,C)=\inf\{r>0,A\subset C+rB\}$ and $B$ is the unit ball.

We assume that a Lyapunov function $V:[0,\,\infty)\times \mathbb{R}^{n}\to \mathbb{R}$ satisfies:
\begin{equation}\label{2.04}
\left\{
\begin{array}{l}
V(\cdot,\cdot) \textup{ is continuous}, \\
\textup{the tube } t \mapsto {\mathcal E}pi V(t,\cdot)=\{(x,v)\in \mathbb{R}^{n}\times \mathbb{R}:\,v\geq V(t,x)\} \textup{ is absolutely continuous.}
\end{array}
\right.
\end{equation}

If $V$ is a locally Lipschitz continuous function then the conditions (\ref{2.04}) are satisfied.

We say that $V$ is positively definite if $V(t,x)>0$ for $x\neq 0$ and $V(t,0)=0$, for $t\geq 0$.\\
We say that $V$ is radially unbounded if there exists an increasing function $p:[0,\infty)\to[0,\infty)$ such that
\begin{equation}\label{2.12}
V(t,x)\geq p(|x|) \textrm{ for every }(t,x)
\end{equation}
and $\lim_{r\to\infty}p(r)=\infty$. This codition has been used in \cite{MP}.

We assume that the function $W:[0,\infty)\times \mathbb{R}^{n}\to \mathbb{R}$ satisfies:
\begin{equation}\label{2.05}
\left\{
\begin{array}{l}
W(t,\cdot) \textup{ is lower semicontinuous for a.a. } t,\\
W(\cdot,\cdot) \textup{ is } {\mathcal L}\times{\mathcal B} \textup{- measurable},\\
|W(t,x)|\leq k(t)(1+|x|) \textup{ for a.a. } t\in[0,\infty) \textup{ and all } x.
\end{array}
\right.
\end{equation}

We assume that the function $g:[0,\infty)\to[0,\infty)$ satisfies
\begin{equation}\label{2.05a}
\left\{
\begin{array}{l}
g(\cdot) \textup{ is continuous and increasing, } g(0)=0 \\
\textup{and } \int_0^b\frac{1}{g(v)}dv<\infty \textup{ for } b>0.
\end{array}
\right.
\end{equation}

Observe that the function $g(v)=v^\alpha$, where $\alpha\in(0,1)$, satisfies (\ref{2.05a}).
The contingent epiderivative of the function $\varphi:\mathbb{R}^n\to \mathbb{R}$ at $x\in \mathbb{R}^n$ in the direction $u\in \mathbb{R}^n$ is defined by
\begin{equation*}
D_{\uparrow}\varphi(x)(u)=\underset{\substack{h \to 0^+ \\ u' \to u}}{\liminf}\frac{\varphi(x+hu')-\varphi(x)}{h}.
\end{equation*}

The contingent hypoderivative of the function $\varphi:\mathbb{R}^n\to \mathbb{R}$ at $x\in \mathbb{R}^n$ in the direction $u\in \mathbb{R}^n$ is defined by
\begin{equation*}
D_{\downarrow}\varphi(x)(u)=\underset{\substack{h \to 0^+ \\ u' \to u}}{\limsup}\frac{\varphi(x+hu')-\varphi(x)}{h}.
\end{equation*}

The properties of contingent epiderivatives and hypoderivatives and its relation to other nonsmooth analysis tools are presented in \cite{AF}.

We say that $V:[0,\infty)\times \mathbb{R}^{n}\to \mathbb{R}$ is a weak Lyapunov function for $F$ with respect to $W$ if there exists a full measure set $D\subset[0,\infty)$ such that
\begin{equation*}
\underset{(t,x) \in D \times \mathbb{R}^n}{\forall} \ \inf_{v \in F(t,x)}D_{\uparrow} V(t,x)(1,v)\leq -W(t,x).
\end{equation*}

Now we recall Theorem 4.4 from \cite{FP}.
\begin{thm}\label{thm2.1}
Assume that $V:[0,\,b)\times \mathbb{R}^{n}\to \mathbb{R}$ is a weak Lyapunov function for $F:[0,b)\times \mathbb{R}^{n}\multimap \mathbb{R}^{n}$ with respect to $W:[0,b)\times \mathbb{R}^{n}\to \mathbb{R}$ and (\ref{2.02}), (\ref{2.04}), (\ref{2.05}) are satisfied, where $b>0$ or $b=\infty$. Then for every initial condition $(t_0,x_0)$ there exists an absolutely continuous solution $x:[t_0,b)\to\mathbb{R}^{n}$ of the Cauchy problem
\begin{equation}\label{2.07}
\left\{\begin{array}{l}
x'(t)\in F(t,x(t)) \\
x(t_0)=x_0
\end{array}
\right.
\end{equation}
such that for every $t_0\leq t<s<b$ we have
\begin{equation}\label{2.08}
V(s,x(s))\leq V(t,x(t))-\int\limits_t^s W(\tau,x(\tau))d\tau.
\end{equation}
\end{thm}
Fact 1 in \cite{Haimo} can be generalized in the following way.
\begin{prop}\label{prop2.2a}
Suppose that a function $g:[0,\infty)\to [0,\infty)$ satisfies (\ref{2.05a}) and a function $c:[0,\infty)\to[0,\infty)$ is locally integrable. If
\begin{equation}\label{c-local}
\int\limits_{t_0}^\infty c(\tau)d\tau > \int\limits_0^{v_0}\frac{1}{g(v)}dv,
\end{equation}
then the solution $\varphi(\cdot)$ of the Cauchy problem
\begin{equation*}
\left\{
\begin{array}{l}
\varphi'(t)=-c(t)g(\varphi(t))\\
\varphi(t_0)=v_0
\end{array}
\right.
\end{equation*}
reaches the origin in the settling time $T(v_0)$ defined by the property $\int_{t_0}^t c(\tau)d\tau<\int_{t_0}^{T(v_0)} c(\tau)d\tau=\int_0^{v_0} \frac{1}{g(v)}dv$ for $t\in[t_0,T(v_0))$.
\end{prop}
\begin{prf}
The absolutely continuous function $G(v)=\int_0^v\frac{1}{g(s)}ds$ satisfies $G(v)>0$ for $v>0$ and $G(0)=\lim_{v\to 0^+}G(v)=0$. If, for some $t>t_0$, $\varphi(s)>0$ for $s\in[t_0,t)$, then
\begin{gather*}
G(\varphi(t))-G(\varphi(t_0))=\int\limits_{t_0}^t \frac{dG(\varphi(s))}{ds}ds=\int\limits_{t_0}^t\frac{1}{g(\varphi(s))}\varphi'(s)ds=\int\limits_{t_0}^t-c(s)ds,
\end{gather*}
which follows $G(\varphi(t))=G(v_0)-\int_{t_0}^t c(s)ds$. Since the function $c$ satisfies (\ref{c-local}), there exists $T>t_0$ such that $\int_{t_0}^T c(s)ds=G(v_0)$. Set $T(v_0)=\inf \left\{T>t_0:\int_{t_0}^T c(s)ds= \right.$ $\left. G(v_0)\vphantom{\frac{1^2}{@^2}}\right\}$. If $s\in[t_0,T(v_0))$, then $\varphi(s)>0$. Suppose to the contrary that $t_1=\inf\{t>t_0:\varphi(t)=0\}<T(v_0)$. Then $\varphi(s)>0$ for $s\in[t_0,t_1)$. Thus $0=G(\varphi(t_1))=G(v_0)-\int_{t_0}^{t_1} c(s)ds$ which is a contradiction.
Thus $G(\varphi(T(v_0)))=0$. Hence, the settling time equals $T(v_0)$.
\end{prf}

If function $c(\cdot)$ satisfies $\int_0^\infty c(\tau)d\tau=\infty$ then (\ref{c-local}) holds true for every $t_0\geq 0$ and $v_0>0$. In Proposition \ref{prop2.2a} it is sufficient to assume that $g(v)>0$ for $v>0$. Then the assumption that $g(\cdot)$ is increasing can be skipped.
\begin{prop}[Comparison Lemma]\label{prop2.2}
Assume that a Carath\'{e}odory function $f:[t_0,\infty)\times[0,\infty)\to [0,\infty)$ is nondecreasing with respect to the second variable, $\varphi_0> 0$ and the Cauchy problem
\begin{equation}\label{2.07a}
\left\{
\begin{array}{l}
\varphi '(t)=-f(t,\varphi(t)) \textrm{ for almost all } t \in [t_0,T] \\
\varphi(t_0)=\varphi_0
\end{array}
\right.
\end{equation}
has an absolutely continuous solution $\varphi:[t_0,T]\to[0,\infty)$. If a bounded by below semicontinuous function $w:[t_0,T]\to[0,\infty)$ satisfies
\begin{equation}\label{2.07b}
w(s)\leq w(t)-\int\limits_t^s f(\tau,w(\tau))d\tau
\end{equation}
for $t_0 \leq t<s\leq T$ and $w(t_0)=\varphi_0$ then
\begin{equation}\label{2.09}
w(t)\leq\varphi(t) \textrm{ for } t\in[t_0,T].
\end{equation}
\end{prop}
\begin{prf}
Suppose that (\ref{2.09}) is not true, i.e. that there exists $t_{2}\in (t_{0},T]$ such that $w(t_{2})>\varphi(t_{2})$. Since the function $w(\cdot)$ is lower semicontinuous and the function $\varphi(\cdot)$ is continuous we obtain that there exists $s \in \lbrack t_{0},t_{2})$ such that
\begin{equation*}
w(t)>\varphi(t) \textrm{ for } t \in (s,t_{2}).
\end{equation*}

Put $t_{1}=\inf \{s \in \lbrack t_{0},t_{2}):w(t)>\varphi(t)$ for all $t \in (s,t_{2}]\}$. If $t_{1}=t_{0}$ then $w(t_{1})=\varphi(t_{1})$. If $t_{0}<t_{1}$ then there exists a sequence $(s_{n})$ convergent to $t_{1}$ from the left such that $w(s_{n}) \leq \varphi(s_{n})$. Thus $w(t_{1}) \leq \liminf_{n\rightarrow \infty}w(s_{n}) \leq \liminf_{n\rightarrow \infty} \varphi(s_{n})=\varphi(t_{1})$. Therefore, for $t \in (t_{1},t_{2})$ we obtain
\begin{equation*}
\varphi(t)<w(t)\leq w(t_{1})-\int\limits_{t_{1}}^{t}f(\tau,w(\tau))d\tau \leq w(t_{1})-\int\limits_{t_{1}}^{t}f(\tau,\varphi(\tau))d\tau \leq \varphi(t_{1})+\int\limits_{t_{1}}^{t}\dot \varphi(\tau)d\tau=\varphi(t),
\end{equation*}
which is a contradiction.
\end{prf}

If an absolutely continuous function $v:[t_0,T]\to [0,\infty)$ satisfies $v'(t)\leq -g(t,v(t))$ for a.a. $t\in[t_0,T]$ then it satisfies (\ref{2.07b}). So, we obtain:
\begin{cor} \label{cor2.3a}
Assume that $f:[t_0,\infty)\times[0,\infty)\to [0,\infty)$ and $\varphi:[t_0,T]\to[0,\infty)$ satisfy the assumptions of Proposition \ref{prop2.2}. If $v:[t_0,T]\to [0,\infty)$ is an absolutely continuous function such that
\begin{equation*}
v'(t)\leq -f(t,v(t)) \textrm{ for almost all } t\in[t_0,T]
\end{equation*}
and $v(t_0)=\varphi(t_0)$, then
\begin{equation*}
v(t) \leq \varphi(t)\textrm{ for } t\in[t_0,T].
\end{equation*}
\end{cor}

The function $\varphi(t)=\left(w(t_{0})^{1-\alpha}-(1-\alpha)\int_{t_{0}}^{t}c(\tau)d\tau \right)^{\frac{1}{1-\alpha}}$ is a solution of the Cauchy problem (\ref{2.07a}) for $f(t,w)=c(t) w^\alpha $ . By Proposition \ref{prop2.2}, we obtain the following corollary.
\begin{cor} \label{cor2.3}
Assume that a lower semicontinuous and bounded function $w:[t_{0},\infty) \rightarrow \lbrack 0,\infty )$ satisfies the integral inequalities
\begin{equation*}
w(s) \leq w(t)-\int\limits_{t}^{s}w(\tau)^{\alpha}c(\tau)d\tau
\end{equation*}
for $t_{0}\leq t<s \leq T$, where $c:[t_{0},\infty) \rightarrow \lbrack 0,\infty)$ is a locally integrable function, $\alpha \in (0,1)$ and
\begin{equation*}
w(t_{0})^{1-\alpha}-(1-\alpha)\int\limits_{t_{0}}^{t}c(\tau)d\tau>0
\end{equation*}
for $t\in(t_{0},T)$.

Then for every $t \in \lbrack t_{0},T]$
\begin{equation} \label{2.11}
w(t)^{1-\alpha} \leq w(t_{0})^{1-\alpha}-(1-\alpha)\int\limits_{t_{0}}^{t}c(\tau)d\tau.
\end{equation}
\end{cor}

Corollary \ref{cor2.3} is a Gronwall type lemma. It is similar to Perov's version of Gronwall lemma (see \cite[Theorem 21]{DR}). The main difference is in the conclusion (\ref{2.11}). For our use we need it in this form (with minus instead of plus).

\begin{thm}[Weak global finite-time stability]\label{thm2.4}
Assume that the right-hand side $F:[0,\infty)\times \mathbb{R}^{n}\multimap \mathbb{R}^{n}$ satisfies (\ref{2.02}), $0\in F(t,0)$ for $t\geq 0$, the function $V:[0,\,\infty)\times \mathbb{R}^{n}\to [0,\infty)$ satisfies (\ref{2.04}), $V(t,0)=0$ for $t\geq 0$ and $V$ is radially unbounded. Suppose that there exists a full measure set $D\subset[0,\infty)$ such that
\begin{gather*}
\underset{(t,x)\in D \times (\mathbb{R}^n \setminus \{0\})}{\forall} \ \inf_{v \in F(t,x)}D_{\uparrow} V(t,x)(1,v) \leq -c(t)g(V(t,x)),
\end{gather*}
where $c:[0,\infty) \rightarrow \lbrack 0,\infty )$ is a locally integrable function such that $\int_0^\infty c(t) dt=\infty$ and $g:[0,\infty)\to[0,\infty)$ satisfies (\ref{2.05a}).

Then for every initial condition $(t_0,x_0)$ there exists a solution $x:[t_0,\infty)\to\mathbb{R}^n$ to (\ref{2.07}) that satisfies $x(T(t_0,x_0))=0$ ,
where the estimation $T(t_0,x_0)$ of the settling time to the solution $x(\cdot)$ is given by
\begin{equation}\label{2.13b}
T(t_0,v_0)=\inf\left \{t>t_0:\;\int\limits_{t_0}^t c(\tau)d\tau \geq \int\limits_0^{V(t_0,x_0)}\frac{1}{g(v)}dv \right\}.
\end{equation}
\end{thm}
\begin{prf}
Fix $(t_0,x_0)$. By (\ref{2.12}), there exists $r>0$ such that the set $\{x:\exists \ t\geq 0, V(t,x)\leq V(t_0,x_0)\}$ is bounded by $r$.

Set $T=T(t_0,x_0)$. We define the function $W(t,x)$ by
\begin{equation*}
W(t,x)=
\left\{
\begin{array}{ll}
c(t) g(V(t,x)) & \textrm{ if } |x|<r+1, \ t \in[t_0,T] \\
0 & \textrm{ elsewhere.}
\end{array}
\right.
\end{equation*}

The function $W$ is Lebesgue-Borel measurable and lower semicontinuous with respect to the state variable. It satisfies the linear growth condition in (\ref{2.05}) for $k(t)=Mc(t)$, where $M=g(\sup\{V(t,x):|x|\leq r+1 \textrm{ and } t\in[t_0,\,T]\})$. Moreover, since $V(t,0)=0$ for $t\geq 0$ and $g(0)=0$ we have $W(t,0)=0$ for $t\geq 0$. Therefore, because $0\in F(t,0)$ for a.a. $t\geq 0$, we easily calculate for these $t$ the following estimation
\begin{equation*}
\inf\limits_{v \in F(t,0)}D_{\uparrow} V(t,0)(1,v) \leq D_{\uparrow} V(t,0)(1,0) \leq \underset{\substack{h \to 0^+}}{\liminf}\frac{V(t+h,0)-V(t,0)}{h}=0=-W(t,0). 
\end{equation*}
Thus $V$ is a weak Lyapunov function for $F$ with respect to $W$.

By Theorem \ref{thm2.1}, there exists a solution $x:[t_0,\infty)\to \mathbb{R}^n$ of (\ref{2.07}) satisfying
(\ref{2.08}).
Since $W\geq 0$ then $V(t,x(t))\leq V(t_0,x_0)$ for $t\geq t_0$. Thus $|x(t)|<r+1$ for $t\geq t_0$ and $W(\tau,x(\tau))=c(\tau)g(V(\tau,\,x(\tau)))$ for $\tau\in[t_0,T]$. So
\begin{equation*}
V(s,x(s))\leq V(t,x(t))-\int\limits_t^s c(\tau)g(V(\tau,x(\tau)))d\tau
\end{equation*}
for $t_0\leq t<s\leq T$.

By Proposition \ref{prop2.2} we have
\begin{equation*}
V(t,x(t))\leq \varphi(t) \textup{ for } t\in[t_0,T],
\end{equation*}
where $\varphi$ is a solution to the Cauchy problem
\begin{equation*}
\left\{
\begin{array}{l}
\varphi'(t)=-c(t)g(\varphi(t))\\
\varphi(t_0)=V(t_0,x_0).
\end{array}
\right.
\end{equation*}

By Proposition \ref{prop2.2a}, we obtain that $\varphi(T)=0$.\\ So, $V(T,x(T))=0$.
\end{prf}

If the function $g(v)$ in Theorem \ref{thm2.4} is of the form $g(v)=v^\alpha$ and $\alpha\in(0,1)$ then
\begin{equation*}
T(t_0,x_0)=\inf\{ t>t_0:\;V(t_0,x_0)^{1-\alpha}=(1-\alpha)\int\limits_{t_0}^t c(\tau)d\tau\}.
\end{equation*}

\begin{thm}[Weak local finite-time stability]\label{thm2.5}
Assume that $\mathcal{O}$ is an open and bounded neighborhood of the origin in $\mathbb{R}^{n}$, the right-hand side $F:[0,\infty)\times \mathcal{O}\multimap \mathbb{R}^{n}$ satisfies (\ref{2.02}), $0\in F(t,0)$ for $t\geq 0$, the bounded function $V:[0,\infty)\times \mathcal{O}\to [0,\infty)$ satisfies (\ref{2.04}) and is positively definite. Suppose that there exists a full measure set $D\subset[0,\infty)$ such that
\begin{equation*}
\underset{(t,x)\in D\times (\mathcal{O} \setminus \{0\})}{\forall} \ \inf\limits_{v\in F(t,x)}D_{\uparrow} V(t,x)(1,v) \leq -c(t)g(V(t,x)),
\end{equation*}
where $c:[0,\infty) \rightarrow \lbrack 0,\infty )$ is a locally integrable function such that for each $t_0\geq 0$ we have $\int_{t_0}^\infty c(t)dt>0$ and $g:[0,\infty)\to[0,\infty)$ satisfies (\ref{2.05a}).

Then for every $t_0\geq 0$ and $\varepsilon>0$ there exists $\delta>0$ such that for every initial condition $|x_0|<\delta$ there exists a solution $x:[t_0,\infty)\to\varepsilon B$ to (\ref{2.07}) that satisfies $x(T(t_0,x_0))=0$, where the estimation $T(t_0,x_0)$ of the settling time is given by (\ref{2.13b}).
\end{thm}
\begin{prf}
Fix $t_0\geq 0$ and put $M_{t_0}=\int_{t_0}^\infty c(t) dt$. By the assumption $M_{t_0}>0$. Since the function $V(t_0,\cdot)$ is continuous and $V(t_0,0)=0$ we can find $R>0$ such that $R\bar B\subset \mathcal{O}$, where $\bar B$ denotes the closed unit ball, and $\rho_0=\sup\{V(t_0,x):|x|\leq R\}$ satisfies inequality $\int_0^{\rho_0}\frac{1}{g(v)}dv< M_{t_0}$.

Define 
\begin{equation*}
\tau(t_0,\rho)=\inf\left\{t>t_0:\int\limits_{t_0}^t c(\tau) d\tau=\int\limits_0^\rho\frac{1}{g(v)}dv\right\}
\end{equation*}
for $\rho\in(0,\rho_0]$.

Let $T_0=\tau(t_0,\rho_0)$. Choose any $\varepsilon \in (0,R)$ and put
\begin{equation*}
V_\varepsilon=\inf\{V(t,x):t\in[t_0, T_0],|x|=\varepsilon\}.
\end{equation*}

We choose $\delta\in (0,\varepsilon)$ such that
\begin{equation*}
\sup\{V(t_0,x):|x|\leq\delta\}<\frac{V_\varepsilon}{2}.
\end{equation*}

Define $\tilde V:[t_0,T_0)\times \mathbb{R}^{n}\to[0,\infty)$, $\tilde F:[t_0,T_0)\times \mathbb{R}^{n}\multimap \mathbb{R}^{n}$ and $W:[t_0,T_0)\times \mathbb{R}^{n}\to \mathbb{R}$ by
\begin{equation*}
\tilde F(t,x)=\left\{
\begin{array}{ll}
F(t,x) & \textrm{ if } V(t,x)<\frac{V_\varepsilon}{2}, \ |x|<\varepsilon, \\
k(t)(1+\varepsilon)\bar B & \textrm{elsewhere},
\end{array}
\right.
\end{equation*}
\begin{equation*}
\tilde V(t,x)=\left\{
\begin{array}{ll}
V(t,x) & \textrm{ if } V(t,x)<\frac{V_\varepsilon}{2}, \ |x|<\varepsilon, \\
\frac{V_\varepsilon}{2} & \textrm{elsewhere},
\end{array}
\right.
\end{equation*}
\begin{equation*}
W(t,x)=\left\{
\begin{array}{ll}
c(t)g(V(t,x)) & \textrm{ if } V(t,x)<\frac{V_\varepsilon}{2}, \ |x|<\varepsilon, \\
0 & \textrm{elsewhere}.
\end{array}
\right.
\end{equation*}

To apply Theorem \ref{thm2.1} we have to show that:
\begin{enumerate}
\item{the function $\tilde V$ is continuous;}
\item{the tube $t\mapsto {\mathcal E}pi \tilde{V}(t,\cdot)$ is absolutely continuous;}
\item{$\tilde V$ is a Lyapunov function for $\tilde F$ with respect to $W$.}
\end{enumerate}

If $|x|=\varepsilon$ then $V(t,x)\geq V_\varepsilon$. Hence the function $\tilde V$ is constant in a neighborhood of $(t,x)$. The continuity of $\tilde V$ in other points is obvious.

To obtain the second property observe that the Hausdorff semi-distance of the epigraphs of the functions $f_C$, $g_C$, where $f_C=\min(f,C)$ and $C$ is a constant, is estimated by the Hausdorff semi-distance of the epigraphs of functions $f,g:\mathcal O\to[0,\infty)$. To verify this statement assume that
\begin{equation}\label{2.13}
\underset{(x_1,y_1)\in {\mathcal E}pi(f)}{\forall} \ \underset{(x_2,y_2)\in {\mathcal E}pi(g)}{\exists} \ \max(|x_2-x_1|,|y_2-y_1|)<\varepsilon.
\end{equation}

Condition (\ref{2.13}) is equivalent to
\begin{equation*}
\underset{x_1\in\mathcal{O}}{\forall} \ \underset{x_2\in\mathcal{O}}{\exists} \ |x_2-x_1|<\varepsilon \textrm{ and } g(x_2)<f(x_1)+\varepsilon.
\end{equation*}

If $f(x_1)\geq C$, then $g_C(x_2)\leq C<f_C(x_1)+\varepsilon$. If $f(x_1)<C$ and $g(x_2)<f(x_1)+\varepsilon$, then $g_C(x_2)\leq g(x_2)<f(x_1)+\varepsilon=f_C(x_1)+\varepsilon$.

Since the domain $\mathcal{O}$ and the function $V$ are bounded, the compact set $K$ in (\ref{2.03}) can be skipped.

To obtain the third condition it is sufficient to show that
\begin{equation*}
\inf\limits_{v\in k(t)(1+\varepsilon)\bar B}D_{\uparrow} \tilde V(t,x)(1,v)\leq -W(t,x)
\end{equation*}
if $V(t,x)=\frac{V_\varepsilon}{2}$, $|x|<\varepsilon$ and $t$ belongs to a full measure set $D$. However, this is obvious, as $\tilde V\leq V$ in a neighborhood of $(t,x)$ and $F(t,x)\subset k(t)(1+\varepsilon)\bar B$.

By Theorem \ref{thm2.1}, for arbitrary initial condition $(t_0,x_0)\in[0,\infty)\times\mathbb{R}^n$ there exists a solution $\bar x:[t_0,T_0)\to \mathbb{R}^n$ to the differential inclusion $x'(t)\in \tilde F(t,x(t))$ that satisfies
\begin{equation*}
\tilde V(t,\bar x(t))+\int\limits_s^t W(\tau,\bar x(\tau))d\tau\leq \tilde V(s,\bar x(s))
\end{equation*}
for every $t_0\leq s<t<T_0$. Since $W\geq 0$ then $\tilde V(t,\bar x(t))\leq \tilde V(t_0,x_0)$.

If $|x_0|<\delta$, then the trajectory $\bar x(t)$ remains in the ball $\varepsilon B$. Indeed, if there exists $t_1\in(t_0,T_0)$ such that $|\bar x(t_1)|\geq \varepsilon$ then there exists $t_2\in(t_0,t_1]$ such that $|x(t_2)|=\varepsilon$. Hence $V(t_2,\bar x(t_2))\geq V_\varepsilon$. So, $\tilde V(t_2,\bar x(t_2))=\frac{V_\varepsilon}{2}$ and further $\frac{V_\varepsilon}{2} \leq \tilde V(t_0,x_0)=V(t_0,x_0)<\frac{V_\varepsilon}{2}$, which is a contradiction. Therefore, the trajectory remains in the area when $\tilde F=F$, $\tilde V=V$ and $W(t,\bar x(t))=-c(t)g(V(t,\bar x(t)))$. We have $\tau(t_0,V(t_0,x_0))\leq T_0$. By Proposition \ref{prop2.2} we obtain that
\begin{equation*}
V(t,\bar x(t))\leq \varphi(t) \mbox{ for $t\in(t_0,\tau(t_0,V(t_0,x_0)))$},
\end{equation*}
where $\varphi$ is a solution to the Cauchy problem
\begin{equation*}
\left\{
\begin{array}{l}
\varphi'(t)=-c(t)g(\varphi(t))\\
\varphi(t_0)=V(t_0,x_0).
\end{array}
\right.
\end{equation*}

By Proposition \ref{prop2.2a}, we obtain $\varphi(\tau(t_0,V(t_0,x_0)))=0$. 

Thus $\bar x(\tau(t_0,V(t_0,x_0)))=0$ and the settling time of the solution $\bar x(\cdot)$ is estimated by $\tau(t_0,V(t_0,x_0))=T(t_0,x_0)$.
\end{prf}
\begin{thm}[Strong local finite-time stability]\label{thm2.5strong}
Assume $\mathcal{O}$ is an open and bounded neighborhood of the origin in $\mathbb{R}^{n}$, the right-hand side $F:[0,\infty)\times \mathcal{O}\multimap \mathbb{R}^{n}$ satisfies (\ref{2.02}), $0\in F(t,0)$ for $t\geq 0$, the bounded function $V:[0,\infty)\times \mathcal{O}\to [0,\infty)$ is locally Lipschitz continuous and positively definite. Suppose that there exists a full measure set $D\subset[0,\infty)$ such that
\begin{equation}
\underset{(t,x)\in D\times (\mathcal{O} \setminus \{0\})}{\forall} \ \sup\limits_{v\in F(t,x)}D_{\downarrow} V(t,x)(1,v) \leq -c(t)g(V(t,x)), \label{DV}
\end{equation}
where $c:[0,\infty) \rightarrow \lbrack 0,\infty )$ is a locally integrable function such that for each $t_0\geq 0$ we have $\int_{t_0}^\infty c(t)dt>0$ and $g:[0,\infty)\to[0,\infty)$ satisfies (\ref{2.05a}).

Then for every $t_0\geq 0$ and $\varepsilon>0$ there exists $\delta>0$ such that for every initial condition $|x_0|<\delta$ every solution $x(\cdot)$ to (\ref{2.07}) is bounded by $\varepsilon$ and satisfies $x(t)=0$ for $t\geq T(t_0,x_0)$, where $T(t_0,x_0)$ is given by (\ref{2.13b}).
\end{thm}
\begin{prf}
Fix a solution $x:[t_0,\,t_1) \rightarrow \mathbb{R}^n$ to (\ref{2.01}).

Since the function $V$ is locally Lipschitz continuous and $x(\cdot)$ is absolutely continuous then the function $v:[t_0,t_1)\to [0,\infty)$ given by $v(t)=V(t,x(t))$ is absolutely continuous. If $t$ belongs to the full measure set $D$ and the derivative $x'(t)=f$ exists and belongs to $F(t,x)$ then the right upper Dini derivative $D^+v(t)=\limsup_{h\to 0 ^+}\frac{v(t+h)-v(t)}{h}$ satisfies
\begin{equation*}
D^+v(t)\leq D_{\downarrow} V(t,x(t))(1,f).
\end{equation*}
Hence, for almost all $t\in[t_0,\infty)$ we have
\begin{equation*}
v'(t)\leq -c(t)g(v(t)).
\end{equation*}
So, the function $v(\cdot)$ is nonincreasing.

By the same construction as in the beginning of the proof of Theorem \ref{thm2.5}, for any fixed $t_0\geq 0$ and $\varepsilon>0$ we can choose $\delta\in (0,\varepsilon)$.
Let $|x_0|<\delta$ and $x(\cdot)$ be a solution to (\ref{2.07}). Using the same argument as in the end of the proof of Theorem \ref{thm2.5} we assure that the solution $x$ is bounded by $\varepsilon$. It means that $x(\cdot)$ is extendable onto $[t_0,\infty)$.

By Corollary \ref{cor2.3a}, we obtain that
\begin{equation*}
V(t, x(t))\leq\varphi(t) \mbox{ for $t\in(t_0,\tau(t_0,V(t_0,x_0)))$},
\end{equation*}
where $\varphi$ is a solution to the Cauchy problem
\begin{equation*}
\left\{
\begin{array}{l}
\varphi'(t)=-c(t)g(\varphi(t))\\
\varphi(t_0)=V(t_0,x_0).
\end{array}
\right.
\end{equation*}

By Proposition \ref{prop2.2a} we obtain $\varphi(\tau(t_0,V(t_0,x_0)))=0$. 

Thus $\bar x(\tau(t_0,V(t_0,x_0)))=0$ and the settling time of $x(\cdot)$ is estimated by $\tau(t_0,V(t_0,x_0))=T(t_0,x_0)$.
\end{prf}
\begin{rem}\
\begin{enumerate}
\item{The regularity of the function $v(\cdot)$ plays very important role in Corollary \ref{cor2.3a}. If the function $v(\cdot)$ is merely of bounded variation then the estimation of their derivative that holds almost everywhere does not allow to estimate its increase (decrease). So, we assume in Theorem \ref{thm2.5strong} that $V$ is locally Lipschitz continuous instead of (\ref{2.04}).}
\item{The contingent conditions describing weak and strong Lyapunov functions can be equivalently formulated with Bouligand tangent cones to the epigraph or hypograph of the function $V$. Namely, the weak Lyapunov function condition is equivalent to the following viability (weak invariance) condition (comp. \cite{AF}, \cite{FP})
\begin{equation*}
(\{1\}\times F(t,x)\times\{-W(t,x)\})\cap T_{\mathcal {E}pi(V)}(t,x,V(t,x)) \neq \varnothing.
\end{equation*}
}
\item{If the right-hand side is single valued, i.e. $F(t,x)=f(t,x)$ and the Cauchy problem
\begin{equation*}
\left \{
\begin{array}{l}
x'(t)=f(t,x(t)) \\
x(t_0)=x_0
\end{array}
\right.
\end{equation*}
has the unique forward solution then the strong and the weak finite-time stability results are equivalent.}
\item{If a right-hand side $f(t,\cdot)$ is discontinuous then we replace the ordinary differential equation $x'=f(t,x)$ by the differential inclusion $x'\in F(t,x)$, where $F(t,\cdot)$ is a Fillipov's regularization of $f(t,\cdot)$. Then the set-valued map $F(t,\cdot)$ is merely upper semicontinuous. In the next section we apply Theorem \ref{thm2.5strong} to obtain finite-time stability of a state discontinuous ordinary differential equation, which describes a Hopfield neural network model.}
\end{enumerate}
\end{rem}
\chapter{Applications to neural networks}
\setcounter{equation}{0}
The finite-time stability have been successfully applied in many fields. The stability of neural networks is a challenging and very important problem in such application as associative memory, signal and image processing, pattern recognition. We need to know whether the neural network we deal with is able to recognize the pattern or to associate data. However, in fact we have to know more: are these recognitions possible in finite-time. The Hopfield neural network is a model network allowing to check whether tools to study it are valuable or not. Thus we investigate the Hopfield neural network of the form (see \cite{H}):
\begin{equation}
c_{i}\frac{dx_{i}(t)}{dt}=\sum\limits_{j=1}^{n}T_{ij}(t)g_{j}(t,x_{j}(t))-\frac{\tilde{h}_{i}(t,x_{i}(t))}{R_{i}}+\tilde{I}_{i}(t), \label{Ho}
\end{equation}
$i=1,\ldots,n$, where $x_{i}\in\mathbb{R}$ is the state of the $i$th neural cell, $g_i(t,x_{i})$ is a nondecreasing activation function of $x_{i}$, $\tilde{h}_i$ measures the rate with which the $i$th unit will reset its potential, $\tilde{I}_{i}$ is an outer input, $T_{ij}$ simulates the connection of the cells, $c_{i}$ is a total input capacitance and $R_{i}$ is an input resistance. We transform (\ref{Ho}) to the following matrix form
\begin{equation}
\dot x(t)=-h(t,x(t))+B(t)g(t,x(t))+I(t), \label{Matrix_nn}
\end{equation}
where $x=(x_{1},\ldots x_{n})$, $h(t,x)=(h_{1}(t,x_{1}),\ldots$, $h_{n}(t,x_{n}))$, $h_{i}=\frac{\tilde{h}_{i}}{R_{i}c_{i}}$, $B(t)=(b_{ij}(t))$, $b_{ij}(t)=\frac{T_{ij}(t)}{c_{i}}$, $g(t,x)=(g_{1}(t,x_{1})$, $
\ldots,g_{n}(t,x_{n}))$, $I=(I_{1},\ldots,I_{n})=\left(\frac{\tilde{I}_{1}}{c_{1}},\ldots,\frac{\tilde{I}_{n}}{c_{n}}\right)$. 

We assume that $B(t)g(t,0)=-I(t)$, $h(t,0)=0$, $h$ is Carath\'{e}odory function, $g_i(t,x_i)$ is measurable in $t$ and nondecreasing with respect to $x_i$, $B$, $I$ are measurable in $t$. The fact that we assume that $g_i(t,\cdot)$ is only nondecreasing means that we admit that it can contain the point of discontinuity. The assumption $B(t)g(t,0)=-I(t)$ means that the bias function $I(t)$ strictly relates (is determined) to the activation function $g$ at the point $(t,0)$ and the weights connection $b_{ij}(t)$ between neurons.

Because we admit that the right-hand side of (\ref{Matrix_nn}) can be discontinuous with respect to $x$, using Filippov's idea (see \cite{F}, \cite{LC}), we consider the following differential inclusion
\begin{equation} \label{Discontinuous_nn}
\dot x(t) \in -h(t,x(t))+B(t) G(t,x(t))+I(t),
\end{equation}
where $G(t,x)=G_1(t,x_1)\times\ldots\times G_n(t,x_n)$ and $G_i(t,x_i)=[g_i(t,x_i^-),g_i(t,x_i^+)]$, where $g_i(t,x_i^-)$ and $g_i(t,x_i^+)$ denote the left-handed limit and right-handed limit (respectively) of the function $g_i(t,\,\cdot)$ at the point $x_i$ (Filippov's regularization of a nondecreasing function is given by one-sided limits).

We shall consider finite-time stability of (\ref{Discontinuous_nn}) in the ball $\mathcal{O} = B(0,\rho)=\left\{x:|x|< \rho \right \}$ for some $\rho \in (0,1)$. Let us assume that there exist integrable functions $a(\cdot)$, $b(\cdot):[0,\infty) \rightarrow [0,\infty)$, $\alpha \in (0,1)$, $\delta>0$ such that
\begin{equation*}
a(t)-b(t)\rho^{2(1-\alpha)}\geq\delta, 
\end{equation*}
\begin{equation}
a(t) |x|^{2\alpha} \leq xh(t,x), \label{ht}
\end{equation}
\begin{equation*}
I(t)=0, 
\end{equation*}
\begin{equation}
xB(t)g(t,x) \leq b(t)|x|^{2} \label{Be}
\end{equation}
for $t \in \lbrack 0,\infty)$ and $x \in B(0,\rho)$. Note that as $\alpha $ belongs to $(0,1)$ therefore $h$ satisfying (\ref{ht}) need not be linear as it is usually assumed, $a(t)$ estimates below the rate with which the $i$th unit will reset its potential. Moreover, in spite that (\ref{Be}) means that $g$ is at most linear but we do not assume that it is Lipschitz. As the function $V$ from (\ref{2.04}) we use the following one:
\begin{equation*}
V(t,x)=C(t)|x| \left(|x|-\exp \left(-|x|^{\alpha -1}\right)\right),
\end{equation*}
for $t \in [0,\infty)$ and $x \in B(0,\rho)$, where $C(t)=\exp(- t)$. Define also $c(t)=\delta C(t)^{1-\alpha}=\delta \exp((\alpha-1)t)$. We see that for each $t_0\geq 0$, $\int_{t_0}^{\infty }c(t)dt>0 $. All assumptions of Theorem \ref{thm2.5strong} will be satisfied if we check condition (\ref{DV}). Since $V$ is smooth in $(0,\infty)\times(0,\setminus \{0\})$, therefore for $v\ \in F(t,x)$ one has
\begin{equation*}
D_{\downarrow} V(t,x)(1,v)=V_{t}(t,x)+V_{x}(t,x)v=V_{t}(t,x)+V_{x}(t,x)(w-h(t,x)),
\end{equation*}
where $w=B(t)f,\ f\in G(t,x)$, i.e.
\begin{equation*}
D_{\downarrow} V(t,x)(1,v)=V_{t}(t,x)+V_{x}(t,x)(B(t)f-h(t,x)),
\end{equation*}
where $f\in G(t,x)$.

To obtain (\ref{DV}) first we show that
\begin{gather}\label{nier-prz}
V_{t}(t,x)+V_{x}(t,x)(B(t)g(t,x)-h(t,x))+c(t)(V(t,x))^{\alpha} \leq 0
\end{gather}
for $t \in [0,\infty)$, $x \in B(0,\rho) \setminus \{0\}$ for sufficiently small $\rho \in (0,1)$. To this aim observe that
\begin{align*}
& V_{t}(t,x)+V_{x}(t,x)(B(t)g(t,x)-h(t,x))+c(t)\left(V(t,x)\right)^{\alpha} \\
&=C'(t)|x| \left(|x|-\exp \left(-|x|^{\alpha-1}\right)\right)+C(t)\left(x \left(2-\frac{1}{|x|}\exp \left(-|x|^{\alpha -1}\right)\right.\right. \\
& \left. \left. -(1-\alpha)|x|^{\alpha-2}\exp \left(-|x|^{\alpha-1}\right)\vphantom{\frac 12}\right) \left(B(t)g(t,x)-h(t,x) \right) \vphantom{\frac 12}\right)+c(t)(C(t))^{\alpha}|x|^\alpha \left(|x|-\exp \left(-|x|^{\alpha -1}\right) \right)^\alpha \\
& \leq -C(t)|x| \left(|x|-\exp \left(-|x|^{\alpha -1}\right)\right)+C(t)|x|^\alpha \left(b(t)|x|^{2-\alpha}-a(t)|x|^\alpha \right) \\
& \left(2-\left(\frac{1}{|x|}+(1-\alpha)|x|^{\alpha -2} \right) \exp \left(-|x|^{\alpha -1} \right) \right)+\delta C(t)|x|^\alpha \left(|x|-\exp \left(-|x|^{\alpha -1}\right) \right)^\alpha \\
& \leq -C(t)|x| \left(|x|-\exp \left(-|x|^{\alpha -1}\right)\right)+C(t)|x|^{2\alpha} \left(b(t)|x|^{2(1-\alpha)}-a(t) \right)+\delta C(t)|x|^{2\alpha} \\
& \leq -C(t)|x| \left(|x|-\exp \left(-|x|^{\alpha -1}\right)\right)-C(t)|x|^{2\alpha}\delta+\delta C(t)|x|^{2\alpha} \leq 0.
\end{align*}

If $f=(f_1,\ldots,f_n)$ is a vertex of rectangular $G(t,x)$ then $f_i\in\{g_i(t,x_i^-)$, $g_i(t,x_i^+)\}$. If $f_i$ equals to the left-handed limit then we choose a sequence $x_{k_i}$ converging to $x_i$ from the left. If $f_i$ equals to the right-handed limit then we choose a sequence $x_{k_i}$ converging to $x_i$ from the right. Set $x_k=(x_{k_1},\ldots,\,x_{k_n})$. Then $\lim_{k\to \infty}g(t,x_k)=f$. The functions $V(t,\cdot)$, $V_t(t,\cdot)$, $V_x(t,\cdot)$ and $h(t,\cdot)$ are continuous. By (\ref{nier-prz}) we have
\begin{equation*}
V_{t}(t,x_k)+V_{x}(t,x_k)(B(t)g(t,x_k)-h(t,x_k))+c(t)(V(t,x_k))^{\alpha} \leq 0.
\end{equation*}

Passing to the limit we obtain that
\begin{equation}\label{nier-prz_f}
V_{t}(t,x)+V_{x}(t,x)(B(t)f-h(t,x))+c(t)(V(t,x))^{\alpha} \leq 0.
\end{equation}

Since the rectangular $G(t,x)$ is a convex hull of its vertices and the function $f\to V_x(t,x)B(t)f$ is a linear functional then (\ref{nier-prz_f}) holds true for $f\in G(t,x)$, which gives (\ref{DV}).

By Theorem \ref{thm2.5strong}, where $g(v)=v^{\alpha}$, we conclude that the Hopfield neural network (\ref{Ho}) is strongly finite-time stable.

\end{document}